\documentclass[12pt]{amsart}

\usepackage[margin=1.15in]{geometry}

\usepackage{amsmath,amscd,amssymb,amsfonts,latexsym,wasysym, mathrsfs, mathtools,hhline,color}
\usepackage[all, cmtip]{xy}
\usepackage{csquotes}
\usepackage{url}

\definecolor{hot}{RGB}{65,105,225}

\usepackage[pagebackref=true,colorlinks=true, linkcolor=hot ,  citecolor=hot, urlcolor=hot]{hyperref}

\usepackage{ textcomp }
\usepackage{ tipa }

\usepackage{graphicx,enumerate}

\theoremstyle{plain}
\newtheorem{theorem}{Theorem}[section]
\newtheorem{prop}[theorem]{Proposition}

\newtheorem{lm}[theorem]{Lemma}

\newtheorem{cor}[theorem]{Corollary}
\newtheorem{conj}[theorem]{Conjecture}

\newtheorem{thrm}[theorem]{Theorem}

\theoremstyle{definition}

\newtheorem{defn}[theorem]{Definition}
\newtheorem{que}[theorem]{Question}
\newtheorem{rmk}[theorem]{Remark}

\newtheorem{ex}[theorem]{Example}
\newtheorem*{ex*}{Example}

\newtheorem*{ack}{Acknowledgements}

\def\be{\begin{equation}}
\def\ee{\end{equation}}

\def\bt{\begin{thrm}}
\def\et{\end{thrm}}

\def\bc{\begin{cor}}
\def\ec{\end{cor}}

\def\br{\begin{rmk}}
\def\er{\end{rmk}}

\def\bp{\begin{prop}}
\def\ep{\end{prop}}

\def\bl{\begin{lm}}
\def\el{\end{lm}}

\def\bex{\begin{ex}}
\def\eex{\end{ex}}

\def\bd{\begin{defn}}
\def\ed{\end{defn}}

\newcommand\sH{{\mathcal H}}

\newcommand\sP{{\mathcal P}}

\newcommand\sF{{\mathcal F}}

\newcommand\sM{{\mathcal M}}

\newcommand\sV{{\mathcal V}}
\newcommand\sL{\mathcal{L}}

\newcommand\kk{{\mathbb{K}}}

\DeclareMathOperator{\codim}{codim}              
                  
\DeclareMathOperator{\rhomo}{RHom}

\DeclareMathOperator{\spec}{Spec}

\DeclareMathOperator{\sFc}{\sF^{\centerdot}}

\DeclareMathOperator{\Alb}{Alb}
\DeclareMathOperator{\alb}{alb}

\DeclareMathOperator{\Hom}{Hom}

\DeclareMathOperator{\im}{Im}
\DeclareMathOperator{\Char}{Char}

\def\bC{\mathbb{C}}

\def\lra{\longrightarrow}

\def\bZ{\mathbb{Z}}

\title[Perverse sheaves on semi-abelian varieties]{Perverse sheaves on semi-abelian varieties -- a survey of properties and applications}

\author{Yongqiang Liu}
\address{BCAM - Basque Center for Applied Mathematics, Mazarredo, 14 E48009 Bilbao, Basque Country - Spain} 
\email{yliu@bcamath.org}
\author{Laurentiu Maxim}
\address{Department of Mathematics,         University of Wisconsin-Madison,  480 Lincoln Drive, Madison WI 53706-1388, USA.}
\email {maxim@math.wisc.edu}
\author{Botong Wang}
\address{Department of Mathematics,         University of Wisconsin-Madison,  480 Lincoln Drive, Madison WI 53706-1388, USA.}
\email {wang@math.wisc.edu}
\date{\today}


\dedicatory{Dedicated to the memory of Prof. \c Stefan Papadima}


\begin{document}

\maketitle
\begin{abstract}  
We survey recent developments in the study of perverse sheaves on semi-abelian varieties. As concrete applications, we discuss various obstructions on the homotopy type of complex algebraic manifolds (expressed in terms of their cohomology jump loci), homological duality properties of complex algebraic manifolds, as well as new topological characterizations of semi-abelian varieties.
\end{abstract}

 \tableofcontents


\section{Introduction}\label{intro}
In recent years, perverse sheaves have become indispensable tools for studying complex algebraic varieties. They have seen important applications in fields like geometry and topology (e.g., the decomposition theorem \cite{BBD} and the topology of complex algebraic maps), representation theory (e.g., proof of the Kazhdan-Lusztig conjecture, proof of the geometrization of the Satake isomorphism, proof of the fundamental lemma in the Langlands program, etc.), or geometric combinatorics (e.g., Stanley's proof of McMullen's conjecture \cite{St}, Huh-Wang's resolution of the Dowling-Wilson ``top-heavy conjecture'' \cite{HW}, etc.); see  
\cite{CM} for a nice survey, and also \cite[Sections 9.3 and 9.4]{Max}. Furthermore, perverse sheaves are the backbone of M. Saito's mixed Hodge module theory \cite{Sa0, Sa1}, a far-reaching generalization of Hodge theory. 

However, despite their fundamental importance, perverse sheaves remain rather misterious objects. In his 1983 ICM address, MacPherson \cite{MP} stated the following:
\begin{displayquote}
{\it ``The category of perverse sheaves is important because of its applications. It would be interesting to understand its structure more directly.''}
\end{displayquote}
Alternative descriptions of perverse sheaves have since been given in various contexts, e.g., in \cite{MV} (in terms of zig-zags), or in \cite{GMV} (using quivers). 
Perverse sheaves on complex affine tori have been studied by Gabber-Loeser \cite{GL} via the Mellin transformation. Perverse sheaves on abelian varieties have been recently described by Schnell in \cite{Sch} via properties of their cohomology jump loci. The references \cite{GL,Sch} are the departure point for much of the  work  in \cite{LMW3, LMW4}, which is surveyed here,  
 where a complete (global) characterization of $\bC$-perverse sheaves on semi-abelian varieties is given in terms of their 
cohomology jump loci. It is also shown  in \cite{LMW4} (and extending our earlier results from \cite{LMW3}) that perverse sheaves on semi-abelian varieties satisfy a set of properties, collectively termed the {\it propagation package}. Besides providing new obstructions on the homotopy type of smooth complex quasi-projective varieties (via their cohomology jump loci), the propagation package also has other important topological applications, such as finiteness properties of Alexander-type invariants, generic vanishing of Novikov homology, or for the study of homological duality properties of complex algebraic varieties (see also \cite{LMW1,LMW2}).

\medskip

The purpose of this survey is to provide a gentle introduction to these  concepts and results, with an emphasis on geometric examples, applications  and open questions.

\begin{ack} We would like to thank the organizers of the conference {\it  Topology and Geometry: A conference in memory of \c Stefan Papadima} (Bucharest, Romania, May 2018), where part of this work was presented. 
Y. Liu was supported by the ERCEA 615655 NMST Consolidator Grant and also by the Basque Government through the BERC 2018-2021 program and by the Spanish Ministry of Science, Innovation and Universities: BCAM Severo Ochoa accreditation SEV-2017-0718. L. Maxim was partially supported by the Simons Foundation Collaboration Grant \#567077 and by CNCS-UEFISCDI grant PN-III-P4-ID-PCE-2016-0030. B. Wang was partially supported by the NSF grant DMS-1701305.
\end{ack}


\section{Preliminaries}\label{prelim}
Our study of perverse sheaves on semi-abelian varieties is partly motivated by the problem  of understanding possible homotopy types of complex algebraic varieties. In this section, we collect the main definitions and technical tools needed to formulate our results from  \cite{LMW3} and \cite{LMW4}.

\subsection{Cohomology jump loci} 
Let $X$ be a smooth connected complex quasi-projective variety with positive first Betti number, i.e., $b_1(X)>0$. The {\it character variety} $\Char(X)$ of $X$ is the identity component of the moduli space of rank-one $\bC$-local systems on $X$, i.e., 
\begin{center}
$\Char(X):= 
\Hom (H_1(X,\bZ)/\text{Torsion}, \bC^*)\cong (\bC^*)^{b_1(X)}.$
\end{center}

\bd The {\it $i$-th cohomology jump locus of $X$} is defined as: 
$${\sV^{i}(X)=\lbrace \rho\in \Char(X) \mid  H^{i}(X,L_{\rho})\neq 0 \rbrace},$$ 
where $L_{\rho}$ is the unique rank-one $\bC$-local system  on $X$ associated to the representation $\rho\in \Char(X) $.\ed

The jump loci $\sV^i(X)$ are closed subvarieties of $\Char(X)$ and homotopy invariants of $X$.  They provide a unifying framework for the study of a host of questions concerning homotopy types of complex algebraic varieties. In particular, they can be used to tackle Serre's problem concerning groups which can be realized as fundamental groups of smooth complex quasi-projective varieties.

\br \label{degree 0} The jump loci $\sV^i(X)$ can in fact be defined for any finite CW complex $X$.  We list a few properties related to the jump loci, which are useful in computations:
\begin{itemize}
\item[(a)] If $X$ is connected, then $\sV^0(X)=\{ \bC_X\}$ consists of just one point, the trivial rank-one local system. 
\item[(b)] If $X$ is a closed oriented manifold of real dimension $n$, then Poincar\'{e} duality yields that $$H^i(X, L_\rho)^\vee\cong H^{n-i}(X, L_{\rho^{-1}})$$
for any $i$.
\item[(c)] For any $\rho \in \Char(X)$ with associated rank-one local system $ L_\rho$, one has: $\chi(X)= \chi(X, L_\rho)$.
\end{itemize}
\er 
We include here the following elementary examples.
\bex \label{circle} If $X=S^1$, Poincar\'{e} duality yields that $\sV^i(X)=\begin{cases}
\{\bC_X\}, &\text{ if } i=  0,1\\
\emptyset, &\text{ otherwise}.
\end{cases}$  
If $X= (S^1)^n$, the K\"{u}nneth Theorem yields that $\sV^i(X)=\begin{cases}
\{\bC_X\}, &\text{ if } 0\leq i\leq n\\
\emptyset, &\text{ otherwise}.
\end{cases}$ 
  
\eex

\bex \label{curve}  Let $X=\Sigma_g$ be a smooth complex projective curve of genus $g\geq 2$.  Then $\chi(X)=2-2g \neq 0$, and Poincar\'{e} duality yields that  $\sV^i(X)=\begin{cases}
\{\bC_X\}, &\text{ if } i=  0, 2,\\
\Char(X), & \text{ if } i=1, \\
\emptyset, &\text{ otherwise}.
\end{cases}$
\eex
 
 \subsection{Semi-abelian varieties}
As it will be discussed below, the classical Albanese map construction allows us to interpret the cohomology jump loci of a smooth connected complex quasi-projective variety as cohomology jump loci of certain constructible complexes of sheaves (or even of perverse sheaves, if the Albanese map is proper) on a semi-abelian variety.
 
Recall that an {\it abelian variety} of dimension $g$ is a compact complex torus $\bC^{g}/{\bZ^{2g}}$, which is also a complex projective variety.
A  {\it semi-abelian variety} $G$ is an abelian complex algebraic group, which is an extension
$$1 \to T \to G \to A \to 1,$$
where $A$ is an abelian variety and $T$ is a complex affine torus. 

(Semi-)abelian varieties are naturally associated to smooth (quasi-)projective varieties via the {\it Albanese map} construction, see \cite{Iit}.
Specifically, if $X$ is a smooth  complex quasi-projective variety, the {\it Albanese map} of $X$ is a morphism $$\alb: X \to \Alb(X)$$ from $X$ to a semi-abelian variety $\Alb(X)$ such that for any morphism $f: X\to G$ to a semi-abelian variety $G$, there exists a unique morphism $g:\Alb(X) \to G$ such that the following diagram commutes: 
\begin{center}
$\xymatrix{
X \ar[rd]^{f}  \ar[r]^{\alb\quad} & \Alb(X) \ar[d]^{g}\\
           &   G 
}$
\end{center} 
Here, $\Alb(X)$ is called the {\it Albanese variety} associated to $X$.

An important property of the Albanese map is that it induces an isomorphism on the torsion-free part of the first integral homology groups, i.e., 
\be\label{eq0}
H_1(X,\bZ)/\text{Torsion} \overset{\cong}\longrightarrow H_1(\Alb(X),\bZ).
\ee
In particular, (\ref{eq0}) yields the identification:
\be\label{eq1}
\Char(X) \cong \Char(\Alb(X)),
\ee
and the equality of Betti numbers $b_1(X)=b_1(\Alb(X))$. Moreover, the dimension of the abelian part of $\Alb(X)$ is $\dim_{\bC} H^{1,0}$ and the dimension of the affine torus part of $\Alb(X)$ is $\dim_{\bC} H^{1,1}$, where $H^{1,0}$ and $H^{1,1}$ denote the weight $(1,0)$ and, resp., weight $(1,1)$ parts of the mixed Hodge structure on $H^1(X)$. Note also that if $X$ is smooth and projective, then $H^1(X)$ carries a pure Hodge structure of weight $1$ and $\Alb(X)$ is an abelian variety. On the other hand, if $X$ is a smooth quasi-projective variety so that  the mixed Hodge structure on $H^1(X)$ is of type $(1,1)$, then $\Alb(X)$ is a complex affine torus. Examples of such smooth varieties with first cohomology group of type $(1,1)$  include complements of essential hyperplane arrangements and complements of   toric hyperplane arrangements.

\subsection{Constructible complexes and perverse sheaves}
{\it Constructible sheaves} are the algebraic counterpart of the decomposition of a variety $X$ into  manifolds pieces. In more detail, a sheaf $\sF$ of $\bC$-vector spaces on $X$ is {\it constructible} if there exists a finite stratification $X=\sqcup_{\lambda} S_{\lambda}$ of $X$ into locally closed smooth subvarieties (called {\it strata}), such that the restriction of $\sF$ to each stratum $S_{\lambda}$ is a $\bC$-local system of finite rank. A complex $\sF^{\centerdot}$ of sheaves of $\bC$-vector spaces on $X$ is called {constructible} if its cohomology sheaves $\sH^i(\sF^{\centerdot})$ are all constructible. 
Denote by $D^b_c(X;\bC)$ the bounded derived category of constructible complexes on $X$, i.e., one identifies complexes containing the same cohomological information. For $\sF^{\centerdot} \in D^b_c(X;\bC)$, only finitely many of its cohomology sheaves are nonzero and, for every $i$, the set ${\it supp} \, \sH^i(\sF^{\centerdot})$, i.e., the closure of the set of points at which the stalk of $\sH^i(\sF^{\centerdot})$ is nonzero, is an algebraic subvariety of $X$. There is also a general form of duality for complexes (called  {\it Verdier duality}) generalizing classical Poincar\'e duality. 

{\it Perverse sheaves} are an important class of constructible complexes, introduced in \cite{BBD} as a formalization of the celebrated Riemann-Hilbert correspondence of Kashiwara \cite{K1,K2} and Mebkhout \cite{Me1, Me2}, which relates the topology of algebraic varieties (intersection homology) and the algebraic theory of differential equations (microlocal calculus and holonomic $D$-modules). 
A {perverse sheaf} on $X$ is a bounded constructible complex $\sP \in D^b_c(X;\bC)$ such that the following vanishing property holds for $\sF^{\centerdot}=\sP$ and for its Verdier dual $\sP^{\vee}$:
$$\dim{\it supp} \, \sH^{-i}(\sF^{\centerdot}) \leq i, \ \ \ {\rm for \ any \ } i \in \bZ.$$
(Unless otherwise specified, throughout the paper $\dim$ will denote the complex dimension of a complex algebraic variety.) 
Perverse sheaves form an abelian category $Perv(X;\bC)$, and can be thought of as singular versions of local systems. Indeed, a local system $L$ on a smooth complex $n$-dimensional variety yields a perverse sheaf $L[n]$; conversely, a perverse sheaf restricts to a (shifted) local system on some dense open subvariety. More details on the theory of perverse sheaves can be found, e.g., in \cite{BBD}, \cite{CM}, \cite{D2}, or \cite{Max}.

We mention here the following important examples of perverse sheaves:
\begin{enumerate}
\item[(i)] If $X$ is a smooth complex algebraic variety of pure dimension $n$, then the shifted constant sheaf ${\bC}_X[n]$ is a perverse sheaf on $X$.
\item[(ii)] If $X$ is a complex algebraic variety of pure dimension $n$, which moreover is a local complete intersection, then ${\bC}_X[n]$ is a perverse sheaf on $X$.
\item[(iii)] If $X$ is a complex algebraic variety of pure dimension $n$, then the Goresky-MacPherson {\it intersection cohomology complex} $IC_X$ (see \cite{GM2}) is a perverse sheaf on $X$. If, moreover, $X$ is a rational homology manifold then $IC_X\simeq {\bC}_X[n]$.
\end{enumerate}

Let $f:Y \to X$ be a proper surjective algebraic morphism, with $Y$ smooth and of (complex) pure dimension $n$. Such a map can be {\it stratified}, which (roughly speaking) means that there is a stratification $X=\sqcup_{\lambda} S_{\lambda}$ of $X$ so that the restriction $f^{-1}(S_{\lambda}) \to S_{\lambda}$ of $f$ above any stratum  $S_{\lambda}$ is a topologically locally trivial fibration. For any $\lambda$, let $s_{\lambda} \in S_{\lambda}$ denote any point in the stratum $S_{\lambda}$. The map $f$ is called {\it semi-small} if, for any $\lambda$, one has:
\be
\dim S_{\lambda} + 2\dim f^{-1}(s_{\lambda}) \leq \dim X.
\ee
Of course, closed embeddings or finite maps are semi-small (onto their image).
The map $f$ is called {\it small} if the above inequality is strict for all positive codimension strata $S_{\lambda}$. Note that a semi-small map is finite on any open stratum of $X$, hence $\dim Y=\dim X$. 
It is also not difficult to see that, in the above notations, the map $f$ is semi-small if, and only if, $Rf_*({\bC}_Y[n])$ is a perverse sheaf on $X$.

\subsection{Jump loci via constructible complexes}
Back to cohomology jump loci, we note that by the projection formula and (\ref{eq1}), for any $\rho \in \Char(X)\cong \Char(\Alb(X))$ we have: 
\be\label{eq2}
H^i(X, L_\rho) \cong H^i(X, \bC_X \otimes L_\rho) 
\cong H^i\big(\Alb(X), (R\alb_* \bC_X)\otimes L_\rho\big).
\ee
If, moreover, the Albanese map $\alb:A \to \Alb(X)$ is proper (e.g., if $X$ is projective), then the  BBDG decomposition theorem \cite{BBD} 
yields that $R \alb_* \bC_X $ is a direct sum of (shifted)  
{perverse sheaves}. In view of (\ref{eq2}), this motivates the study of {\it cohomology jump loci of constructible complexes} (resp., {\it perverse sheaves}) {\it  on semi-abelian varieties}.

\bd Let $\sFc \in D^b_c(G;\bC)$ be a bounded constructible complex of $\bC$-sheaves on a semi-abelian variety $G$. 
The {\it degree $i$ cohomology jump locus of $\sFc$} is defined as:
$$\sV^i(G,\sF): = \{\rho \in \Char(G)\mid H^i(G, \sFc \otimes_\bC L_\rho) \neq 0\}.$$
\ed

With the above definition, (\ref{eq2}) yields that:
\be\label{eq3}
{\sV^i(X)= \sV^i(\Alb(X), R\alb_* \bC_X)}.
\ee



\section{Mellin transformation and applications}
\subsection{Mellin transformation}
Let $G$ be a semi-abelian variety defined by an extension
$$1 \to T \to G \to A \to 1,$$
where $A$ is an abelian variety of dimension $g$ and $T\cong(\bC^*)^m$ is a complex affine torus of dimension $m$. Set 
$$\Gamma_G:=\bC[\pi_1(G)]\cong \bC[t_1^{\pm 1}, \cdots, t_{m+2g}^{\pm 1}],$$
and note that $$\Char(G)= \spec \Gamma_G.$$
Let $\sL_G$ be the (universal) rank-one local system of $\Gamma_G$-modules on $G$, defined by mapping the generators of $\pi_1(G)\cong \bZ^{m+2g}$ to the multiplication by the corresponding variables of  $\Gamma_G$.

\bd\cite{GL} The {\it Mellin transformation}  $\sM_{\ast}: D^{b}_{c}(G; \bC) \to D^{b}_{coh}(\Gamma_G)$ is given by
$${\sM_\ast(\sFc) := Ra_\ast( \sL_G\otimes_{\bC}\sFc)}, $$
where $a :G\to pt$ is the constant map to a point space, and $D^{b}_{coh}(\Gamma_G)$ denotes  the bounded coherent complexes of $\Gamma_G$-modules. 
\ed

The Mellin transformation can be used to completely characterize perverse sheaves on complex affine tori. More precisely, one has the following result due to Gabber-Loeser \cite[Theorem 3.4.1 and Theorem 3.4.7]{GL} in the $\ell$-adic context, and then extended to the present form in  \cite[Theorem 3.2]{LMW3}: 
\bt\label{GL}
 A constructible complex $\sFc \in D^b_c(T;\bC)$ on a complex affine torus $T$  is perverse if, and only if, \begin{center} $H^i(\sM_\ast(\sFc))=0$ for all $i \neq 0$.\end{center} 
\et

In the context of abelian varieties, the Mellin transformation was considered in \cite{BSS} for proving generic vanishing statements:
\bt\cite{BSS}  For any $\bC$-perverse sheaf $\sP$ on an  {abelian variety} $A$, one has:\begin{center} $H^i( \sM_* (\sP))=0$ for $i<0$.\end{center} and \begin{center}
$H^i (D_{\Gamma_A} (\sM_*(\sP))) =0$ for $i<0$.
\end{center}
Here $D_{\Gamma_A}(-):=\rhomo_{\Gamma_A} (-,\Gamma_A) $ is the dualizing functor for the ring $\Gamma_A=\bC[\pi_1(A)]$.
\et

This result was recently extented to the semi-abelian context in \cite{LMW4}, by induction on the dimension of the complex affine torus $T$:
\bt\cite[Theorem 4.3]{LMW4} \label{Mellin}  For any $\bC$-perverse sheaf $\sP$  on a  {semi-abelian variety} $G$, one has: \begin{center} $H^i( \sM_* (\sP))=0$ for $i<0$, \end{center} 
and \begin{center}
$H^i (D_{\Gamma_G} (\sM_*(\sP))) =0$ for $i<0$.
\end{center}
Here $D_{\Gamma_G}(-):=\rhomo_{\Gamma_G} (-,\Gamma_G) $ is the dualizing functor for $\Gamma_G$-modules.
\et

\subsection{Propagation and codimension lower bound}
The Mellin transformation was used in \cite{LMW3}, \cite{LMW4} to translate the question of understanding the cohomology jump loci of a constructible complex to a problem in commutative algebra.
Specifically, by the projection formula, cohomology jump loci of $\sFc \in D^b_c(G;\bC)$  are determined by those of the Mellin transformation $\sM_\ast(\sFc)$ of $\sFc$ as follows:
\be\label{eq5}{\sV^i(G,\sFc)=\sV^i(\sM_\ast(\sFc))}.\ee
Here, if $R$ is a Noetherian domain and $E^{\centerdot}$ is a bounded complex of $R$-modules with finitely generated cohomology, we set
$$ \sV^{i} (E^\centerdot):= \{ \chi \in \spec R \mid H^i (F^\centerdot \otimes _R R/\chi) \neq 0 \},$$ 
 with $F^\centerdot$ a bounded above finitely generated {\it free} resolution of $E^\centerdot$. We also have the following:

\bp \label{exact} \cite[Lemma 2.8]{BSS} Let $R$ be a regular Noetherian domain and $E^{\centerdot}$ a bounded complex of $R$-modules with finitely generated cohomology. Then $H^i(E^{\centerdot} )=0$ for $i<0$ if, and only if, $\codim \sV^{ -i}(E^{\centerdot}) \geq i$ for $i\geq 0$.
\ep
 
By using the identification (\ref{eq5}), together with  Proposition \ref{exact} and  standard facts from commutative algebra, we proved in \cite{LMW4} the following result as a direct consequence of Theorem \ref{Mellin}:
\bt\cite[Theorem 4.6]{LMW4}\label{pro}     For any $\bC$-perverse sheaf $\sP$ on a semi-abelian variety $G$, the cohomology jump loci of $\sP$ satisfy the following properties:
\begin{enumerate}
\item[(i)] {\it Propagation}:
$$
{ \sV^{-m-g}(G,\sP) \subseteq \cdots   \subseteq \sV^{0}(G,\sP) \supseteq \sV^{1}(G,\sP) \supseteq \cdots \supseteq \sV^{g}(G,\sP)}.$$
Moreover, $\sV^i(G,\sP)= \emptyset$ if $i\notin [-m-g, g]$.
\item[(ii)] {\it Codimension lower bound}: for any $i\geq 0$, $$ {\codim \sV^{ i}(G,\sP) \geq i} \ \
{\rm and} \  \ {\codim \sV^{ -i}(G,\sP) \geq i}.$$
\end{enumerate}
\et

\br \label{coefficient} Cohomology jump loci of constructible complexes can be defined for any field coefficients, and in fact Theorem \ref{Mellin} and Theorem \ref{pro} hold in this generality, see \cite[Theorem 4.3, Theorem 4.6]{LMW4}. 
\er

The result of Theorem \ref{pro} can be viewed as a topological counterpart of similar properties satisfied by the Green-Lazarsfeld algebraic jump loci of topologically trivial line bundles, see \cite{GLa, GrLa}.

\br[Equivalent formulation of the propagation property]
Let $\sP$ be a $\bC$-perverse sheaf so that not all  $H^j(G, \sP)$ are zero.  Let 
\begin{center} $k_+:=\max\{j\mid H^j(G, \sP)\neq 0\} \text{ and } k_-:=\min\{j \mid H^j(G, \sP)\neq 0\}.$\end{center} 
The propagation property (i) is equivalent to: $k_+\geq 0$, $k_-\leq 0$ and
$$H^j(G, \sP)\neq 0 \ \iff \ k_-\leq j\leq k_+. $$
Moreover,  if $\sP$ is a semi-simple perverse sheaf on an abelian variety $G=A$, then the Hard Lefschetz Theorem  yields that $k_-=-k_+$. 
Then by the relative Hard Lefschetz theorem for the Albanese map, a similar interval of nonvanishing cohomology can be obtained for a simple perverse sheaf on any smooth complex projective variety; compare with \cite[Corollary 1]{We16b}. 
\er


\subsection{Application: generic vanishing and signed Euler characteristic}
An important consequence of Theorem \ref{pro} is the following generic vanishing statement:
\bc\label{c1} For any $\bC$-perverse sheaf $\sP$ on a semi-abelian variety $G$,  $$H^{i}(G, \sP\otimes_{\bC} L_\rho)=0$$ for any generic rank-one $\bC$-local system $L_\rho$ and all $i\neq 0$. 
In particular,
 \be\label{eq6}  \chi(G,\sP):=\sum_i (-1)^i \dim_{\bC} H^i(G,\sP)\geq 0.\ee
Moreover, $ \chi(G,\sP)=0$ if, and only if, $\sV^0(G,\sP) \neq \Char(G)$.
\ec
\begin{proof}
If $\rho \notin \bigcup_{i\neq 0} \sV^i(G,\sP)$, then $H^{i}(G, \sP\otimes_{\bC} L_\rho)=0$ for $i\neq 0$. The codimension bound in Theorem \ref{pro} shows that $\bigcup_{i\neq 0} \sV^i(G,\sP)$ has at least codimension one in $\Char(G)$, so the generic vanishing result follows. Moreover, for such a generic choice of character $\rho$, we have:
 $$\chi(G,\sP)=\chi(G,\sP\otimes_{\bC} L_\rho)=\dim_{\bC} H^0(G, \sP\otimes_{\bC} L_\rho )\geq 0.$$ 
The last claim follows directly from the propagation property.
\end{proof}

The above generic vanishing statement was originally proved by other methods in 
 \cite[Theorem 2.1]{Kra} in the $\ell$-adic context and further generalized to arbitrary field coefficients in \cite[Theorem 1.1]{LMW2}. For abelian varieties, generic vanishing results were obtained in \cite[Theorem 1.1]{KW}, \cite[Corollary 7.5]{Sch}, \cite[Vanishing Theorem]{We16} or \cite[Theorem 1.1]{BSS}. The signed Euler characteristic property (\ref{eq6}) is originally due to Franecki and Kapranov \cite[Corollary 1.4]{FK}.  
 
\br The properties of perverse sheaves from Theorem \ref{pro} and Corollary \ref{c1} are collectively termed the {\it propagation package} for perverse sheaves on semi-abelian varieties.\er
 
 \bex\label{exs} Let $Z$ be a closed pure-dimensional subvariety of a semi-abelian variety $G$, with inclusion map $j:Z \hookrightarrow G$. Then $j_*IC_Z \in Perv(G;\bC)$, so (\ref{eq6}) yields that 
 \be\label{eq6ih}
 (-1)^{\dim Z} \chi_{IH}(Z) \geq 0,
 \ee
 where $\chi_{IH}(Z):=\sum_{i \geq 0} (-1)^i \dim_{\bC} IH^i(Z,\bC)$ is the intersection homology Euler characteristic of $Z$. If $G=A$ is an abelian variety, a geometric proof of (\ref{eq6ih}) was recently obtained in \cite{EGM} by using the stratified Morse theory of Georesky-MacPherson \cite{GMs}; see also \cite{BW15} for the case when $G=T$ is a complex affine torus.
 
 If the closed subvariety $Z\subseteq G$ is smooth, or a rational homology manifold, or a local complete intersection, then $j_*\bC_Z[\dim Z] \in Perv(G;\bC)$, so (\ref{eq6}) yields in these cases that 
 \be\label{eq6h}
 (-1)^{\dim Z} \chi(Z) \geq 0.
 \ee
 For a geometric proof of (\ref{eq6h}) in the smooth case see \cite{LMW1}, where the non-proper Morse theory of Palais-Smale \cite{PS} is used; see also \cite{Huh} for the case when $G=T$ is a complex affine torus. If $G=A$ is an abelian variety, then (\ref{eq6h}) is proved in the local complete intersection case in \cite{EGM} by stratified Morse theory techniques.
 \eex
 
 \medskip
 
 For completeness, and also in relation with Example \ref{exs}, we also mention here the following version of the signed Euler characteristic property, as a special case of positivity results for characteristic classes obtained recently in \cite{AMSS}, and which is already implicit in \cite{FK} and  \cite{TS}:
 \bt\cite{AMSS,FK,TS}\label{amss}
 If $Z$ is a closed subvariety of an abelian variety $A$, and $\varphi$ is a constructible function on $Z$ with effective characteristic cycle, then 
  \be\label{eq6g}
\chi(Z,\varphi) \geq 0.
 \ee
 \et
In order to explain the meaning of Theorem \ref{amss} and how it relates to (\ref{eq6ih}) and (\ref{eq6h}), let us first recall some terminology (where we can replace $A$ by any smooth complex algebraic variety).
A function $\varphi:Z \to \bZ$ is {\it constructible} if there is a Whitney stratification $\mathcal{Z}$ of $Z$ so that $\varphi$ is constant on any stratum $S \in \mathcal{Z}$. The Euler characteristic of such a constructible function $\varphi$ is defined by
$$\chi(Z,\varphi):=\sum_{S \in \mathcal{Z}} \chi(S) \cdot \varphi(S),$$
where $\varphi(S)$ denotes the (constant) value of $\varphi$ on the stratum $S$.
 
 To any constructible complex $\sFc \in D^b_c(Z;\bC)$ one associates a constructible function $\chi_{st}(\sFc)$ defined by taking stalkwise Euler characteristics, i.e., $$\chi_{st}(\sFc)(x):=\chi(\sFc_x).$$ For example, $\chi_{st}(\bC_Z)=1_Z$ and we let $ic_Z:=\chi_{st}(IC_Z)$. Note that if $\varphi=\chi_{st}(\sFc)$, then $\chi(Z,\varphi)=\chi(Z,\sFc)$.  Another important example of a constructible function is the {\it local Euler obstruction} function $Eu_Z$ of MacPherson \cite{MP0}, which is an essential ingredient in the definition of Chern classes for singular varieties.
 
 To any constructible function $\varphi$ on $Z \subseteq A$, one associates a {\it characteristic cycle} \be\label{ccd}CC(\varphi)=\sum_{Y \subseteq Z} a_Y [T^*_YA] \in L(T^*A|_Z),\ee for uniquely determined closed subvarieties $Y$ of $Z$ and nonzero integer coefficients $a_Y$. Here, $L(T^*A|_Z)$ denotes the group of conic Langrangean cycles in the cotangent bundle $T^*A$ with support on $Z$, and $T^*_YA$ is the conormal space of $Y$, i.e., the closure in $T^*A$ of the conormal bundle of the smooth locus of $Y$.  For example, if $Z \subseteq A$ is a closed and irreducible subvariety, then:
 \be\label{cc}
 CC(Eu_Z)=(-1)^{\dim Z}[T^*_ZA].
 \ee
 In particular, by (\ref{cc}), the coefficients $a_Y$ in (\ref{ccd}) are determined by:
 \be\label{coe}
 0 \neq \varphi =\sum_Y a_Y (-1)^{\dim Y}Eu_Y.
 \ee
 The characteristic cycle of $\varphi\neq 0$ is {\it effective} if and only if the coefficients $a_Y$ are all positive. For example, if $\varphi=(-1)^{\dim Z}Eu_Z$, then (\ref{cc}) yields that $CC(\varphi)=[T^*_ZA]$, which is clearly effective. Therefore, Theorem \ref{amss} yields the following:
 \bc
 If $Z$ is a closed irreducible subvariety of an abelian variety $A$, then
 \be
 (-1)^{\dim Z}\chi(Z,Eu_Z) \geq 0.
 \ee
 \ec
  It is also known that if $0\neq \varphi=\chi_{st}(\sP)$ is a nontrivial constructible function associated to a peverse sheaf $\sP$ on $Z$, then $CC(\varphi)$ is effective. In particular, $CC(ic_Z)$ is effective, so Theorem \ref{amss} recovers (\ref{eq6ih}) in the abelian context. Similarly, if $Z\subseteq A$ is smooth, or a rational homology manifold, or a local complete intersection, then $\bC_Z[\dim Z]\in Perv(Z)$, so Theorem \ref{amss} recovers (\ref{eq6h}) in the abelian context.
  
  Finally, it is shown in \cite[Proposition 1.14]{AMSS} that if $p:G \to A$ is the projection from a semi-abelian variety onto its abelian part, then $p$ has the following property: if a constructible function $\varphi $ on $G$ is associated to a perverse sheaf, then its pushforward $p_*(\varphi)$ (as a constructible function on $A$) also comes from a perverse sheaf, so $p_*(\varphi)$ has an effective characteristic cycle. In particular, an application of Theorem \ref{amss} yields that if $\sP$ is a perverse sheaf on $G$ with associated constructible function $\varphi$ then:
  $$\chi(G,\sP)=\chi(G,\varphi)=\chi(A,p_*(\varphi)) \geq 0,$$
  thus reproving (\ref{eq6}); see \cite[Corollary 1.16]{AMSS}.




\subsection{Application: homological duality}
Abelian duality spaces were introduced in \cite{DSY}, by analogy with the {duality spaces} of Bieri-Eckmann \cite{BE}. They are especially useful for explaining some previously conjectural behavior of the cohomology of abelian representations.

Let $X$ be a finite connected CW complex, with $\pi:=\pi_1(X)$ and $H:=H_1(X,\bZ)$.  
There is a canonical $\bZ [H]$-local system on $X$, whose monodromy action is given by the composition of the quotient $\pi \to \pi^{ab}=H$ with the natural multiplication $H\times \bZ[H] \to \bZ[H]$. 
\bd\cite{DSY}\label{abdsp} A finite connected CW complex $X$ is an {\it abelian duality space} of dimension $n$ if $H^i(X, \bZ[H])=0$ for all $i\neq n$, and $H^n(X, \bZ[H])$ is non-zero and torsion-free (as a $\bZ$-module).
\ed
When $X$ is compact, Definition \ref{abdsp} is equivalent to the requirement that the integral compactly supported cohomology of the universal abelian cover of an abelian duality space $X$ is concentrated in a single dimension, $n$, where it is torsion-free.

A {\it duality space} \cite{BE} is defined similarly, by using $\pi$ instead of $H$ in the above definition (hence by replacing the universal abelian cover by the universal cover). 

Duality spaces and abelian duality spaces enjoy (homological) duality properties similar to the Serre duality for projective varieties.

\medskip


A very natural question to investigate is the following  {\it realization problem} for (abelian) duality spaces:
\begin{que}\label{q2}
Which smooth complex quasi-projective varieties are (abelian) duality spaces?
\end{que}

\noindent For example, there are plenty of smooth complex projective varieties which are duality spaces (e.g., all smooth complex projective curves of genus $g \geq 1$), but, as it will be explained below, 
smooth complex projective varieties which are abelian duality spaces are rather rare (e.g., among smooth complex projective curves, only genus $1$ curves are abelian duality spaces). 

Obstructions to the realizability of abelian duality spaces referred to in  Question \ref{q2} are provided by the following result (see \cite{DSY}, \cite{LMW3}):
\bt\label{ad}
Let $X$ be an abelian duality space of dimension $n$. For simplicity, assume that $H_1(X,\bZ)$ is free\footnote{The statement remains true without any assumption on $H_1(X,\bZ)$, but in this generality one needs to consider all characters of $\pi_1(X)$, i.e., all elements of $\Hom(\pi_1(X), \bC^*)$, and the corresponding cohomology jump loci; 
see \cite[Theorem 4.6]{LMW3} for details.}. Then the cohomology jump loci of $X$ satisfy the propagation package, i.e., the following properties hold:
\begin{enumerate}
\item[(i)] {\it Propagation property}: $$\sV^{n}(X) \supseteq \sV^{n-1}(X) \supseteq \cdots \supseteq \sV^{0}(X).$$
\item[(ii)] {\it Codimension lower bound}: for any $i\geq 0$, $$ \codim \sV^{n-i}(X)= b_1(X)-\dim  \sV^{n-i}(X) \geq i.$$
\item[(iii)]  {\it Surviving property}: If  $V$ is an irreducible component of $\sV^n(X)$ of codimension $d$, then $V \subset \sV^{n-d}(X).$ In particular, $V$ is also an irreducible component of $\sV^{n-d}(X)$.
\item[(iv)] {\it Generic vanishing}: $H^{i}(X,  L_\rho)=0$ for $\rho \in \Char(X)$ generic and all $i\neq n$ .
\item[(v)] {\it Signed Euler characteristic property}:  $$  (-1)^n\chi(X)\geq 0.$$
\item[(vi)] {\it Betti property}: $b_i(X) > 0$, for $0 \leq i \leq n$, and $b_1(X) \geq  n$.
\end{enumerate} 
 \et
\begin{proof}
We only indicate here a proof for the Betti property (vi).
By Remark \ref{degree 0}(a), we have that $\sV^0(X)=\{ \bC_X\}$ is just one point. Then the propagation property (ii) yields that $\bC_X  \in \sV^i(X)$ for $0\leq i \leq n$, hence $b_i(X)>0$ for $0\leq i \leq n$. Furthermore, the codimension bound (ii) for $\sV^0(X)$ shows that $b_1(X)-\dim \sV^0(X) \geq n$, hence $b_1(X)\geq n$.
\end{proof}
As a direct application of the propagation package of Theorem \ref{ad},  in \cite{LMW3} we gave a complete answer to the abelian duality realization problem of Question \ref{q2} in the projective setting:
\bt\cite[Theorem 6.1]{LMW3} \label{abp}
A smooth complex projective variety $X$ is an abelian duality space if, and only if, $X$ is an abelian variety. Moreover, a compact K\"ahler manifold is an abelian duality space if, and only if, it is a compact complex torus.
\et
\noindent Indeed, in this case, Poincar\'e duality imposes additional obstructions on the cohomology jump loci of $X$, 
which, when combined with the propagation property, imply that the Albanese map $\alb:X\to \Alb(X)$ is a finite covering, hence an isomorphism (compare also with Question \ref{q1} below).

Theorem \ref{abp}  provides a new topological characterization of compact complex tori and, respectively, abelian varieties in terms of homological duality properties, and can be viewed as a topological counterpart of the Chen-Hacon {\it birational} characterization of abelian varieties \cite{CH}.

In relation to Theorem \ref{abp}, one can also ask the following:
\begin{que}\label{q4} Does there exist a closed orientable manifold that is an abelian duality space, but not a real torus?\end{que}

Furthermore, by interpreting the abelian duality property in terms of the Mellin transformation, we identifed in \cite{LMW4} a very general class of complex algebraic varieties which can be realized as abelian duality spaces:
 \bt\cite[Theorem 7.11]{LMW4} \label{abel}
Let $X$ be an $n$-dimensional smooth complex quasi-projective variety, which is homotopy equivalent to a CW complex of real dimension $n$ (e.g., $X$ is affine).  Assume that  the Albanese map $\alb$ of $X$ is proper and semi-small (e.g., a finite map or a closed embedding), or $\alb$ is quasi-finite. Then $X$ is an abelian duality space of dimension $n$. 
\et

\begin{proof} Let us briefly sketch the main idea of the proof.
For simplicity, assume that $H_1(X,\bZ)$ is a free abelian group. The assumption that $X$ is homotopy equivalent to an $n$-dimensional CW complex implies that $H^i(X, \bZ[H])=0$ for all $i>n$. 
On the other hand,  it is easy to show that $$H^i(\sM_*(R \alb_* \kk_X))\cong H^i(X, \kk[H]) $$ for any $i\in \bZ$. 
Here $\kk$ is any field of coefficients, and $\kk_X$ is the corresponding constant sheaf on $X$. Since the hypotheses of the theorem imply that $R \alb_* (\kk_X[n])$ is $\kk$-perverse on $\Alb(X)$, Theorem \ref{Mellin} and Remark \ref{coefficient} yield that $H^i(X, \kk[H])=0$ for $i<n$. Using a universal coefficient theorem type of argument, one can then show that  $H^i(X, \bZ[H])=0$ for $i\neq n$ and $H^n(X, \bZ[H])$ is torsion free as a $\bZ$-module.
\end{proof}

Let us next discuss a few examples. 
\bex[Very affine manifolds] \label{affine} Let $X$ be an $n$-dimensional {\it very affine manifold}, i.e., a smooth closed subvariety of a complex affine torus $T=(\bC^*)^m$ (e.g., the complement of an essential hyperplane arrangement or the complement of a  toric arrangement).  
The closed embedding  $i:X\hookrightarrow T$ is a proper semi-small map, and hence the Albanese map $\alb:X \to \Alb(X)$ is also proper and semi-small.  Since $X$ is affine and $n$-dimensional, it has the homotopy type of a finite CW complex of real dimension $n$. Therefore, we get by Theorem \ref{abel} that $X$ is an abelian duality space of dimension $n$. In particular, very affine manifolds satisfy the propagation package of Theorem \ref{ad}.\eex 

\bex \label{ample}  Let $Y$ be a smooth complex projective variety, and let $\sL$ be a very ample line bundle on $Y$. Consider an $N$-dimensional sub-linear system $\vert E\vert $ of $\vert \sL \vert$ such that $E$ is base point free over $Y$. Then a basis $\{s_0, s_1, \cdots, s_N\}$ of $E$ gives  a well-defined morphism  $$\varphi_{\vert E \vert}:   Y \to \mathbb{CP}^N .$$ Each  $\{ s_i=0 \} $ defines a  hypersurface $V_i$ in $Y$. In particular, $\bigcap_{i=0}^N V_i = \emptyset$. Then $\varphi_{\vert E \vert}$ is a finite morphism; for a proof, see \cite[Example 5.4]{LMW3}.
Taking the restriction $f$ of $\varphi_{\vert E \vert}$ over  $X= Y\setminus \bigcup_{i=0}^N V_i $, we get a map
\begin{align*}
 f:=\left(\dfrac{s_1}{s_0}, \dfrac{s_2}{s_0}, \cdots,\dfrac{s_N}{s_0}\right): X 
 \lra T=(\bC^*)^N \end{align*}
which is finite, hence  proper and semi-small. This implies that the Albanese map $\alb$ is also proper and semi-small.
Theorem \ref{abel} yields that $X$ is an abelian duality space. 

This example shows that there are plenty of smooth quasi-projective varieties which are abelian duality space. If all hypersurfaces $V_i$ are smooth and intersect locally like hyperplanes, then this example is a special case of \cite[Theorem 1.1]{DS}. However, here we do not make any assumptions on the singularities of the $V_i$'s and of their intersections.
\eex

\bex[Elliptic arrangement complements]  \label{Elliptic}
It was shown in \cite{DS,DSY} that the complement of an {\it elliptic arrangement} is an abelian duality space. This fact also follows from Theorem \ref{abel} as follows. Let $E$ be an elliptic curve, and let $\mathcal{A}$ be an essential elliptic arrangement in $E^n$ with complement $X:=E^n \setminus \mathcal{A}$. Then $X$ is a complex $n$-dimensional affine variety, and by the universal property of the Albanese map, the natural embedding $X \hookrightarrow E^n$ factorizes through $\alb:X \to \Alb(X)$. Hence the Albanese map $\alb:X \to \Alb(X)$ is also an embedding (hence, in particular, quasi-finite). So Theorem \ref{abel} applies to show that $X$ is an abelian duality space of dimension $n$.
\eex

The affine assumption is a sufficient but not a necessary condition for an $n$-dimensional smooth complex quasi-projective variety to be homotopy equivalent to a finite CW- complex of real dimension $n$. The following is an example of an $n$-dimensional smooth complex quasi-projective variety which is an abelian duality space of dimension $n$, but which is not affine.
\bex
Let $X$ be the blowup of $(\bC^*)^2$ at a point.  The Albanese map for $X$ is the blowdown map $X\to (\bC^*)^2$, which is proper and semi-small.  Moreover, $X$ is homotopy equivalent to the $2$-dimensional CW-complex $T^2 \vee S^2$, where $T^2=S^1\times S^1$ is the real $2$-dimensional torus.  Thus, Theorem \ref{abel} shows that $X$ is an abelian duality space of dimension $2$. However, $X$ is not affine because it contains a closed subvariety $\mathbb{C}\mathbb{P}^1$.
\eex

If an $n$-dimensional smooth complex quasi-projective variety $X$ is an abelian duality space and $\alb: X\to \Alb(X)$ is proper, then $\alb$ is semi-small (onto its image). Indeed, by using the decomposition theorem and the relative hard Lefschetz theorem, one can readily see that $R\alb_* (\bC_X[n])$ is a perverse sheaf.  We conjecture that, in some sense, the converse of Theorem \ref{abel} is also true. 
\begin{conj} Let $X$ be an $n$-dimensional smooth complex quasi-projective variety with proper Albanese map $\alb : X \to \Alb(X)$. Then $X$ is an abelian duality space of dimension $n$ if, and only if,  $\alb$ is semi-small and $X $ is homotopy equivalent to a finite  CW complex of real dimension $n$.  
\end{conj}

It is also natural  to investigate topological properties of the Bieri-Eckmann duality spaces. We formulate here the following:
\begin{que}\label{q3}
If $X$ is a smooth complex quasi-projective variety which is a duality space, what can be said about its topology? For example, is the Euler characteristic of $X$ signed? 
\end{que}

\noindent 
The second part of Question \ref{q3} is motivated by the {\it Hopf and Singer conjectures}, e.g., see \cite[Ch.11]{Lu} for an overview. Hopf conjectured that if $M$ is a closed manifold of real dimension $2n$, 
with negative sectional curvature, then $(-1)^n \chi(M) \geq 0$.
Jost-Zuo \cite{JZ} proved Hopf's conjecture in the K\"ahler context. More precisely, 
they showed that if $M$ is 
a compact K\"ahler manifold of complex dimension $n$ and non-positive 
sectional curvature, then $(-1)^n \chi(M) \geq 0$. 
Hopf's conjecture was strengthened by Singer, who claimed that if $M^{2n}$ 
is a closed {\it aspherical} (i.e., it has a contractible universal cover) manifold, then $(-1)^n \chi(M) \geq 0$. 
(Note that if $M$ carries a Riemannian metric with non-positive sectional 
curvature, then $M$ is aspherical by Hadamard's Theorem.)  
In \cite{LMW3}, we proved the following result:
\bp\label{p4} A closed oriented manifold is a duality space if and only if it is aspherical.
\ep
\noindent In particular, Proposition \ref{p4} provides a reformulation of Singer's conjecture from the point of view of homological duality.

\section{Characterization of perverse sheaves and applications}
\subsection{Characterization of perverse sheaves }

In this last section, we touch upon MacPherson's quote mentioned in the Introduction, and present a complete (global) characterization of perverse sheaves on semi-abelian varieties; see \cite{LMW4} for complete details.  
Motivation is also provided by the following result of Schnell:
\bt\cite{Sch} \label{Sch}
If $A$ is an abelian variety and $\sFc\in D^b_c(A;\bC)$, then  
$\sFc$ is perverse if, and only if,  for any $i \in \bZ$, 
$\codim \sV^{i}(A,\sFc) \geq \vert 2i \vert.$
\et

Furthermore, as a consequence of Theorem \ref{GL}, Proposition \ref{exact} and Artin's vanishing theorem, one gets the following:
\bt \cite[Corollary 6.8]{LMW4}\label{LMWaf}
$\sFc\in D^b_c(T; \bC)$ is perverse on a complex affine torus $T$ if, and only if, 
\begin{enumerate}
\item[(i)] For any $i>0$: $\sV^i(T, \sFc)=\emptyset$, and 
\item[(ii)] For any $i\leq 0$: $\codim \sV^i(T, \sFc)\geq -i$. 
\end{enumerate}
\et

In order to unify and generalize the results of Theorems \ref{Sch} and \ref{LMWaf} to the semi-abelian context, we need to recall a structure result 
for the cohomology jump loci of constructible complexes on a semi-abelian variety $G$, which  has been recently obtained by Budur-Wang in \cite{BW17}; see also \cite[Theorem 2.2]{Sch} for the abelian case.
\bd\label{lin} A closed irreducible subvariety $V$ of $\Char(G)$ is called {\it linear} if
there exists a short exact sequence of semi-abelian varieties 
\begin{center}
$1 \to G''(V)\to G\overset{q}{\to} G'(V) \to 1$
\end{center} and  some $\rho\in \Char(G)$ such that
\begin{center}
$V:=\rho \cdot \im (q^\#: \Char(G'(V)) \to \Char(G))$.
\end{center} 
Here $G''(V)$ and $G'(V)$ depend on $V$, and $q^\#$ is corresponding map induced by $q:G\to G'(V)$.
\ed
With the above definition, one has the following result:
\bt\cite{BW17} \label{bw}
For any $\sFc \in D^b_c(G;\bC)$, each jump locus $\sV^i(G,\sFc)$ is a finite union of  linear subvarieties of $\Char(G)$.
\et

\bd Let $G$ be a semi-abelian variety and let $V$ be an irreducible linear subvariety of $\Char(G)$. In the notations of Definion \ref{lin}, let  
$T''(V)$ and $A''(V)$ denote the complex affine torus and, resp., the  abelian variety part of $G''(V)$. Define:
\begin{itemize}
\item[] {\it abelian codimension}: $\codim_a V:=\dim A''(V),$ 
\item[] {\it semi-abelian codimension}: $\codim_{sa} V:=\dim G''(V).$
\end{itemize}
These definitions generalize to reducible subvarieties by taking the minimum among all irreducible components. 

\ed

\br Let $V$ be a nonempty linear subvariety of $\spec \Gamma_G$. 
\begin{enumerate}
\item If $G=T$ is a complex affine torus, then: $\codim_{sa} V =\codim V$ and $\codim_{a} V=0$. 
\item If $G=A$ is an abelian variety, we have: $\codim_{sa}(V)=\codim_{a}(V)=\frac{1}{2}\codim(V)$. 
In particular, Theorem \ref{Sch} can be restated as follows: $\sFc\in D^b_c(A;\bC)$ is perverse if, and only if,  for any $i \in \bZ$, 
$\codim_a \sV^{i}(A,\sFc) \geq \vert i \vert.$
\end{enumerate}
\er

In the above notations, the following generalization of Schnell's result was obtained in \cite{LMW4}:
\bt\cite[Theorem 6.6]{LMW4} \label{main}
A constructible complex $\sFc\in D^b_c(G; \bC)$ is perverse on $G$ if, and only if, 
\begin{enumerate}
\item[(i)] $\codim_{a} \sV^i(G, \sFc) \geq i$ for any $i\geq 0$, and
\item[(ii)] $\codim_{sa} \sV^i(G, \sFc) \geq -i$ for any $i\leq 0$. 
\end{enumerate}
\et

\begin{proof} 
The ``only if'' part is proved by induction on $\dim T$, using Theorem \ref{Sch} as the beginning step of the induction process. For the ``if'' part, one shows that the two codimension lower bounds in the statement of Theorem \ref{main} are sharp. For complete details, see \cite[Section 6]{LMW4}.
\end{proof}

\subsection{Application: cohomology jump loci of quasi-projective manifolds}

The results of Theorem \ref{pro}, Theorem \ref{main} and Corollary \ref{c1} can be directly applied for the study of cohomology jump loci $\sV^i(X) \subseteq \Char(X)=\Char(\Alb(X))$ of a smooth quasi-projective variety $X$. More precisely, we have the following
 (see \cite[Corollary 7.2]{LMW4} for a more general formulation):

\bc
Let  $X$ be a  smooth quasi-projective variety of complex dimension $n$.  Assume that 
$R\alb_*( \bC_X[n])$ is a perverse sheaf on $\Alb(X)$ (e.g., $\alb$ is proper and semi-small, or $\alb$ is quasi-finite).
Then the cohomology jump loci $\sV^i(X) $ have the following properties:
\begin{itemize}
\item[(1)] {\it Propagation property}:
$$ 
  \sV^{n}(X) \supseteq \sV^{n-1}(X)\supseteq \cdots \supseteq \sV^{0}(X) = \{ 1 \};
$$
$$ 
  \sV^{n}(X) \supseteq \sV^{n+1}(X) \supseteq \cdots \supseteq \sV^{2n}(X)  .
$$
\item[(2)] {\it Codimension lower bound}: for any $ i\geq 0$, 
$$ \codim_{sa} \sV^{n-i}(X) \geq i \ \ {\rm and } \ \ \codim_a \sV^{n+i}(X) \geq i.$$
\item[(3)] {\it Generic vanishing:} $H^{i}(X,  L_\rho)=0$ 
for generic $\rho \in \Char(X)$ and all $ i \neq n$.
\item[(4)] {\it Signed Euler characteristic property}:  $(-1)^n  \chi(X) \geq 0$.
\item[(5)] {\it Betti property}: $b_i(X)>0$ for any $i\in [0, n]$, and $b_1(X)\geq n$.
\end{itemize} 
\ec


 \bex Let $X$ be a smooth closed sub-variety of a semi-abelian variety $G$ (e.g., $X$ is a very affine manifold as in Example \ref{affine}, or a smooth projective curve as in Example \ref{curve}). The closed embedding  $i:X\hookrightarrow G$ is a proper semi-small map, and hence the Albanese map $\alb:X \to \Alb(X)$ is also proper and semi-small. Then $R \alb_* (\bC_X[\dim X])$ is a perverse sheaf on $\Alb(X)$. Hence the jump loci of $X$ satisfy the properties listed in the above corollary.  \eex 
 
 
\subsection{Application: topological characterization of semi-abelian varieties}
The Structure Theorem \ref{bw} and the propagation package of Theorem \ref{pro} and Corollary \ref{c1} can be used to give the following topological characterization of semi-abelian varieties (see \cite[Proposition 7.7]{LMW4}):
\bc\label{c3.5} \label{homotopy} Let $X$ be a smooth quasi-projective variety with proper Albanese map (e.g., $X$ is projective), and assume that $X$ is homotopy equivalent to a real torus. Then $X$ is isomorphic to a semi-abelian variety. 
\ec
\begin{proof} Assume $X$ has complex dimension $n$.
By the decomposition theorem \cite{BBD}, $R\alb_* (\bC_X[n])$ is a direct sum of shifted semi-simple perverse sheaves on $\Alb(X)$. Denote by $\mathcal{S}$ the collection of all simple summands appearing (up to a shift) in $R\alb_* (\bC_X[n])$. Then, by using Theorem \ref{pro}(i) and the identification (\ref{eq3}), we have that \be\label{un} \bigcup_{i=0}^{2n} \sV^i(X) = \bigcup_{\sP\in \mathcal{S}} \sV^0(\Alb(X),\sP).\ee  
Since $X$ is homotopy equivalent to a real torus, Example \ref{circle} shows that  $\bigcup_{i=0}^{2n} \sV^i(X)$ is just an isolated point. Hence, by (\ref{un}), for every simple perverse sheaf $\sP\in \mathcal{S}$, the jump locus $\sV^0(\Alb(X),\sP)$ is exactly this isolated point, so Corollary \ref{c1} yields that $\chi(\Alb(X),\sP)=0$. 

Simple perverse sheaves with zero Euler characteristic on semi-abelian varieties are completely described in \cite[Theorem 5.5]{LMW4} by using the structure Theorem \ref{bw} and the propagation package for their cohomology jump loci. In particular, it follows that for any $\sP \in \mathcal{S}$ we get:  $\sP=\bC_{\Alb(X)}[\dim \Alb(X)]$.  So  $R\alb_* \bC_X$ is a direct sum of shifted rank-one constant sheaves on $\Alb(X)$. Since $X$ and $\Alb(X)$ are both homotopy equivalent to tori, and since $b_1(X)=b_1(\Alb(X))$, we have that $b_i(X)=b_i(\Alb(X))$ for any $i$. Therefore, $$R\alb_* \bC_X\cong \bC_{\Alb(X)}.$$ Since $\alb$ is proper, it follows that all fibers of $\alb$ are zero-dimensional. Then it can be seen easily that $\alb$ is in fact an isomorphism.
\end{proof}

Under the assumptions of Corollary \ref{homotopy}, we showed that the Albanese map 
 $\alb: X\to \Alb(X)$ is an isomorphism. 
This motivates the following:
\begin{que}\label{q1}
Find necessary and sufficient conditions on the cohomology jump loci of $X$ so that the Albanese map $\alb:X \to \Alb(X)$ is an isomorphism,  or a dominant map, or a submersion. 
\end{que}
\noindent For example, we showed in \cite[Corollary 7.4]{LMW4} that if $X$ is an $n$-dimensional smooth quasi-projective variety so that $\bigcup_{i=0}^{2n} \sV^i(X)$ has an isolated point, then $\alb:X \to \Alb(X)$ is a dominant map. Furthermore, we showed in 
\cite[Corollary 7.6]{LMW4} that if the Albanese map of $X$ is proper and semi-small, then in order to conclude that $\alb$ is an isomorphism it is sufficient (and necessary) to ask that $\sV^n(X)$ consists of finitely many points.

\medskip

Before stating the next and final question, we recall the following:
\bd 
A proper morphism of complex spaces $f : X \to Y$ is called a {\it homotopy fiber bundle} if $Y$ has an open cover $Y = \bigcup_j U_j$ such that for every $j$ and for every $y \in U_j$ the inclusion
$$f^{-1}(y) \to f^{-1}(U_j) $$ is a homotopy equivalence.
\ed 
The following is a simplified version of a question of Bobadilla and Koll\'{a}r \cite[Question 4.2 and Conjecture 3]{BK}:
\begin{que}
Let $X$ be a smooth complex projective variety. Assume that the universal free abelian cover of $X$ is homotopy equivalent to a finite CW complex. Does this imply that $\alb: X \to \Alb(X)$ is a homotopy fiber bundle, or even a differentiable fiber bundle?
\end{que}
While this question is still open in general,  Corollary \ref{homotopy} can be viewed as a partial answer.


\end{document}